\documentclass[12pt]{article}
\def\dateref{}
\usepackage{amssymb, amsmath, amsthm, amscd, graphicx, dsfont}

\newtheorem{pro}{Proposition}
\newtheorem{thm}{Theorem}
\newtheorem{defn}{Definition}

\def\cE{{\mathcal E}}
\def\cP{{\mathcal P}}
\def\Z{{\bf Z}}

\def\colon{\,{:}\;}
\def\prf{\medbreak\noindent{\bf Proof.}\enspace}
\def\qed{\hspace*{\fill}\hbox{\vrule height 7pt \kern-.3pt
     \vbox{\hrule width 7pt
     \kern6.6pt\hrule width 7pt }\kern-.3pt\vrule height 7pt
     }\par}
\def\ra{\rightarrow}

\begin{document}

\title{\bf Weak Gibbs measures and large deviations}

\author{C.-E. Pfister\footnote{E-mail: charles.pfister@epfl.ch}\\
Section of Mathematics,
Faculty of Basic Sciences, EPFL\\
CH-1015 Lausanne, Switzerland
\and W.G. Sullivan\footnote{E-mail:
wgs@maths.ucd.ie}\\
               Department of Mathematics, UCD,\\
               Belfield, Dublin 4, Ireland}

 \date{\dateref}

\maketitle

\noindent
{\bf Abstract:\,}
Let $(X,T)$ be a dynamical system, where $X$ is a compact metric space and $T\colon X\ra X$ a continuous onto map.
For weak Gibbs measures we prove large deviations estimates.

\section{Introduction}
\setcounter{equation}{0}

In \cite{PS1} a general method for proving large deviations estimates for dynamical systems (X,T) is developed. In this  note we make the connection with the main results of \cite{PS1} and the notion of weak Gibbs measures, which was not explicit in the original paper.

Let $X$ be a compact metric space and $T\colon X\ra X$ a continuous map which is onto. $M_1(X)$ is the set of Borel probability measures on $X$ (with weak convergence topology) and $M_1(X,T)$ the subset of $T$-invariant probability measures.
Let $x\in X$ and
$$
\cE_n(x):=\frac{1}{n}\sum_{k=0}^{n-1}\delta_{T^kx}\,.
$$
The metric entropy of $\nu\in M_1(X,T)$ is denoted $h(T,\nu)$ and
$B_m(x,\varepsilon)$ is the dynamical ball
$\{y\in X\colon d(T^kx,T^ky)\leq\varepsilon\,,\;k=0,\ldots,m-1\}$.
There are several variants in the literature for the definition of weak Gibbs measures
(see e.g. \cite{BV} and \cite{Yu}).
In this paper a weak Gibbs measure is defined as follows.

\begin{defn}\label{multi-defn2}
Let $\varphi\in C(X)$.
A  probability measure $\nu$ is a \emph{weak Gibbs measure for $\varphi$} if $\forall \delta>0$  $\exists \varepsilon_\delta>0$ such that for
 $0<\varepsilon\leq\varepsilon_\delta$ $\exists N_{\delta,\varepsilon}<\infty$, $\forall m\geq N_{\delta,\varepsilon}$, $\forall x\in X$,
$$
-\delta\leq \frac{1}{m}\ln\nu(B_m(x,\varepsilon))-\int \varphi\,d\cE_m(x)\leq \delta\,.
$$
\end{defn}

The set of weak Gibbs measures for a given $\varphi$ is convex (possibly empty).
Gibbs measures as defined in \cite{Bo} (see \cite{Bo}, theorem 1.2) and quasi-Gibbs measures (see \cite{HR}, proposition  2.1) are examples of weak Gibbs measures since these measures satisfy the stronger inequalities: there exists
$0<\varepsilon_0<\infty$ such that for
$0<\varepsilon\leq\varepsilon_0$ $\exists K_{\varepsilon}<\infty$,
 $\forall m$, $\forall x\in X$,
$$
-\frac{K_{\varepsilon}}{m}\leq \frac{1}{m}\ln\nu(B_m(x,\varepsilon))-\int \varphi\,d\cE_m(x)\leq\frac{K_{\varepsilon}}{m}\,.
$$

\section{Results}
\setcounter{equation}{0}

If $\nu$ is a weak Gibbs measure, then
\begin{eqnarray*}
0&=
\lim_{\varepsilon\downarrow 0}\liminf_m\inf_{x\in X}\Big(\frac{1}{m}\ln\nu(B_m(x,\varepsilon))-\int \varphi\,d\cE_m(x)\Big)\\
&=
\lim_{\varepsilon\downarrow 0}\limsup_m\sup_{x\in X}\Big(\frac{1}{m}\ln\nu(B_m(x,\varepsilon))-\int \varphi\,d\cE_m(x)\Big)\,,
\end{eqnarray*}
that is, $-\varphi$ is a lower, respectively upper, energy function for $\nu$ in the sense of \cite{PS1} 
(definitions 3.2 and 3.4). Indeed, in \cite{PS1} a function $e$ on $X$ is called a \emph{lower energy function for} $\nu$ if it is
upper semi-continuous and
\begin{equation}\label{eq1}
\lim_{\varepsilon\downarrow 0}\liminf_m\inf_{x\in X}\Big(\frac{1}{m}\ln\nu(B_m(x,\varepsilon))+\int e\,d\cE_m(x)\Big)\geq 0\,.
\end{equation}
It is called an \emph{upper energy function for} $\nu$ if it is lower semi-continuous, bounded and
\begin{equation}\label{eq2}
\lim_{\varepsilon\downarrow 0}\limsup_m\sup_{x\in X}\Big(\frac{1}{m}\ln\nu(B_m(x,\varepsilon))+\int e\,d\cE_m(x)\Big)\leq 0 \,.
\end{equation}
The terminology used in \cite{PS1} comes from statistical mechanics.

\begin{pro}\label{pro1}
If the continuous function $e$ verifies (\ref{eq1}) and (\ref{eq2}), then $\nu$ is a weak Gibbs measure for 
$\varphi=-e$.
\end{pro}

\prf
For any $\delta>0$, if $\varepsilon$ is small enough and $m$ large enough,
\begin{eqnarray*}
-\delta  &\leq \inf_{m\geq N_{\delta,\varepsilon}}
\inf_{x\in X}\Big(\frac{1}{m}\ln\nu(B_m(x,\varepsilon))-\int \varphi\,d\cE_m(x)\Big)\\
&\leq
\liminf_m\Big(\frac{1}{m}\ln\nu(B_m(x,\varepsilon))-\int \varphi\,d\cE_m(x)\Big)\\
&\leq
\limsup_m \Big(\frac{1}{m}\ln\nu(B_m(x,\varepsilon))-\int \varphi\,d\cE_m(x)\Big)\\
&\leq \sup_{m\geq N_{\delta,\varepsilon}}
\sup_{x\in X}\Big(\frac{1}{m}\ln\nu(B_m(x,\varepsilon))-\int \varphi\,d\cE_m(x)\Big)\leq \delta\,,
\end{eqnarray*}
so that for $\forall m\geq N_{\delta,\varepsilon}$ and $\forall x\in X$
$$
-\delta\leq \frac{1}{m}\ln\nu(B_m(x,\varepsilon))-\int \varphi\,d\cE_m(x)\leq \delta\,.
$$
\qed

For any dynamical system $(X,T)$ and any weak Gibbs measure the following large deviations estimates are true.

\begin{thm}\label{thm2}
Let $\nu$ be a weak Gibbs measure for $\varphi$.

1.  If $G\subset M_1(X)$ is open,  then for any ergodic probability measure $\rho\in G$
$$
\liminf_m\frac{1}{m}\ln\nu(\cE_m\in G)\geq h(T,\rho)+\int \varphi\,d\rho\,.
$$

2. If  $F\subset M_1(X)$ is convex and closed, then
$$
\limsup_m\frac{1}{m}\ln\nu(\cE_m\in F)
\leq \sup_{\rho\in F\cap M_1(X,T)}\big(h(T,\rho)+\int \varphi\,d\rho\big)\,.
$$
\end{thm}

\prf
Proposition 3.1 and theorem 3.2 in \cite{PS1}. \qed

\begin{pro}\label{pro2}
If $\nu$ is a weak Gibbs measure for $\varphi$, then the topological pressure $P(\varphi)=0$.
\end{pro}

\prf
This is an immediate consequence from theorem \ref{thm2}, theorem 9.10 and corollary 9.10.1 in \cite{Wa}.
Let $G=F=M(X)$. Then
$$
P(\varphi)=\sup_{\rho\,{\rm ergodic}}\big(h(T,\rho)+
\int \varphi\,d\rho\big )
\leq 0 \leq
\sup_{\rho\in M_1(X,T)}\big(h(T,\rho)+
\int \varphi\,d\rho\big )=P(\varphi)\,.
$$
\qed

The following hypothesis about the entropy-map $h(T,\cdot)$ and the dynamical system $(X,T)$
are sufficient to obtain a full large deviations principle.

\begin{thm}
Let $\nu$ be a weak Gibbs measure for $\varphi$.
If  the entropy map $h(T,\cdot)$ is upper semi-continuous, then
for $F\subset M_1(X)$ closed
$$
\limsup_m\frac{1}{m}\ln\nu(\cE_m\in F)
\leq \!\!\!\!\sup_{\rho\in F\cap M_1(X,T)}\big(h(T,\rho)+\int \varphi\,d\rho\big)\,.
$$
If the ergodic measures are entropy dense, then for $G\subset M_1(X)$ open
$$
\liminf_m\frac{1}{m}\ln\nu(\cE_m\in G)
\geq \!\!\!\!\sup_{\rho\in G\cap M_1(X,T)}\big(h(T,\rho)+
\int \varphi\,d\rho\big)\,.
$$
\end{thm}

\prf Theorems 3.1 and 3.2 in \cite{PS1}. \qed

Entropy density of
 the ergodic measures means  (\cite{PS1}): for any $\mu\in M_1(X,T)$, any neighbourhood $N$ of $\mu$ and any $h^*<h(T,\mu)$, there exists an ergodic measure $\rho\in N$ such that
$h(T,\rho)\geq h^*$. Entropy density is true under various types of specifications properties for the dynamical system $(X,T)$, see e.g. \cite{PS1}, \cite{PS2},  \cite{CTY}, \cite{KLO} and \cite{GK}.

\begin{pro}\label{pro4}
If $\nu\in M_1(X,T)$ is a weak Gibbs measure for $\varphi$, then it is an equilibrium measure for $\varphi$.
\end{pro}

\prf
By definition an equilibrium measure $\mu\in M_1(X,T)$ for a continuous function $f$
satisfies  the variational principle
$$
P(f)=\sup\Big\{h(T,\rho)+\int f\,d\rho\colon \rho\in M_1(X,T)\Big\}
=h(T,\mu)+\int f\,d\mu\,.
$$
Since $P(\varphi)=0$,  $h(T,\nu)\leq -\int \varphi\,d\nu$. Since $\nu$ is a weak Gibbs measure for $\varphi$,
\begin{eqnarray*}
\limsup_m\int \varphi\,d\cE_m(x)&=\lim_{\varepsilon\downarrow 0}\limsup_m
\frac{1}{m}\ln\nu(B_m(x,\varepsilon))\\
\liminf_m\int \varphi\,d\cE_m(x)&=\lim_{\varepsilon\downarrow 0}\liminf_m
\frac{1}{m}\ln\nu(B_m(x,\varepsilon))\,.
\end{eqnarray*}
By the ergodic theorem there exists an integrable function $\varphi^*$ such that
$$
\lim_m \int \varphi\,d\cE_m(x)=\varphi^*(x)\quad\nu-{\rm a.s.}
$$
and
$$
\int \varphi\,d\nu=\int \varphi^*\,d\nu\,.
$$
Therefore
$$
h(T,\nu)\leq -\int \varphi(x)\,d\nu(x)=\int\big(-\lim_{\varepsilon\downarrow 0}\limsup_m
\frac{1}{m}\ln\nu(B_m(x,\varepsilon))\big)\,d\nu(x)\,.
$$

Let $\cP=\{A_1,\ldots,A_p\}$ be a finite measurable partition of $X$,
$\max_i{\rm diam}A_i<\varepsilon$. For $x\in X$, let $\cP^n(x)$ be the element of the partition $\cP^n=\cP\vee T^{-1}\cP\vee \cdots\vee T^{-n+1}\cP$ containing $x$. By the
McMillan-Breiman theorem
$$
h_\cP(x):=\lim_n-\frac{1}{n}\ln \nu(\cP^n(x))\quad \nu-{\rm a.s.}
$$
and
$$
\int h_\cP(x)\,d\nu(x)=h_\cP(T,\nu)\,,
$$
where
$$
h_\cP(T,\nu)=\lim_n\Big(-\frac{1}{n}\sum_{B\in\cP^n}\nu(B)\ln \nu(B)\Big)\leq h(T,\nu)\,.
$$
Since $B_n(x,\varepsilon)\supset \cP^n(x)$, for any $\varepsilon>0$,
$$
\int \big(-\limsup_m \frac{1}{m}\ln\nu(B_m(x,\varepsilon))\big)\,d\nu(x)
\leq \int h_\cP(x)\,d\nu(x)\leq h(T,\nu)\,,
$$
so that -$\int \varphi\,d\nu\leq h(T,\nu)$.
\qed

\bigskip
\noindent
{\bf Concluding remark\,}
The results in \cite{PS1} are proven  for  continuous $\Z_+^d$-actions or $\Z^d$-actions on $X$.
The  results of this note are also true for these cases. The empirical measure $\cE_n(x)$ and the dynamical ball $B_n(x,\varepsilon)$ are defined as in \cite{PS1}.

\end{document}